\documentclass[a4paper,12pt,reqno]{amsart}
\usepackage[utf8]{inputenc}
\usepackage{amssymb}
\usepackage{amsmath}
\usepackage{makecell}
\usepackage{multirow}
\usepackage{hhline}
\usepackage{amsthm}
\usepackage{amsfonts}
\usepackage{mathtools}
\usepackage{todonotes}
\usepackage{comment}
\usepackage{mathrsfs}
\usepackage{enumitem}
\usepackage[all,cmtip]{xy}
\usepackage{ctable}
\usepackage[hyphens]{url}
\usepackage{hyperref}
\usepackage[hyphenbreaks]{breakurl}
\usepackage{cleveref}
\usepackage{lipsum}
\usepackage{subcaption}
\usepackage{tabularx}

\usepackage[british]{babel}
\usepackage[margin=1.2in]{geometry}
\usepackage{eqparbox}

\usepackage{pbox}

\newcommand{\field}[1]{\mathbb{#1}}  
\newcommand{\Q}{\field{Q}} 

\DeclareMathOperator{\rk}{rk}

\DeclareMathOperator{\Reg}{\textup{\textsf{Reg}}}

\DeclareFontFamily{U}{wncy}{}
\DeclareFontShape{U}{wncy}{m}{n}{<->wncyr10}{}
\DeclareSymbolFont{mcy}{U}{wncy}{m}{n}
\DeclareMathSymbol{\Sha}{\mathord}{mcy}{"58}

\newtheorem{lemma}{Lemma}

\theoremstyle{definition}

\newtheorem{experiment}[lemma]{Experiment}

\numberwithin{lemma}{section}
\numberwithin{equation}{section} 
\numberwithin{figure}{section}

\title[Machine-Learning Sha]{Machine Learning Approaches to the Shafarevich-Tate Group of Elliptic Curves}
\author[A. Babei]{Angelica Babei}
\address{Angelica Babei \\
Cambridge, MA 02139\\
USA}
\email{babeiangelica@gmail.com}

\author[B.S. Banwait]{Barinder S. Banwait}
\address{Barinder S. Banwait \\
Department of Mathematics \& Statistics\\  
Boston University\\
665 Commonwealth Ave\\
Boston, MA 02215\\
USA; \& \\
Department of Mathematics\\  
Massachusetts Institute of Technology\\
77 Massachusetts Avenue\\
Cambridge, MA 02139\\
USA (visitor)}
\email{barinder@mit.edu}

\author[A. Fong]{AJ Fong}
\address{AJ Fong \\
Department of Pure Mathematics \\
University of Waterloo \\
Waterloo, Ontario \\
N2L 3G1 \\
Canada}
\email{aj.fong@uwaterloo.ca}

\author[X. Huang]{Xiaoyu Huang}
\address{Xiaoyu Huang\\
Department of Mathematics \\ 
Temple University\\
Wachman Hall\\
Philadelphia, PA 19122\\
USA
}
\email{xiaoyu.huang@temple.edu}

\author[D. Singh]{Deependra Singh}
\address{Deependra Singh \\
Department of Mathematics \\
Emory University \\
Atlanta, GA 30322 \\
USA}
\email{deependra.singh@emory.edu}

\date{}

\makeatletter
\providecommand\@dotsep{5}
\renewcommand{\listoftodos}[1][\@todonotes@todolistname]{%
  \@starttoc{tdo}{#1}}
\makeatother

\subjclass[2010]
{11G05  (primary), 
14Q05,   
68T05.   
(secondary)
}

\begin{document}

\maketitle

\begin{abstract}
We train machine learning models to predict the order of the Shafa\-revich-Tate group $\Sha$ of an elliptic curve over $\Q$. Building on earlier work of He, Lee, and Oliver, we show that a feed-forward neural network classifier trained on subsets of the invariants arising in the Birch--Swinnerton-Dyer conjectural formula yields higher accuracies ($> 0.9$) than any model previously studied. In addition, we develop a regression model that may be used to predict orders of $\Sha$ not seen during training and apply this to the elliptic curve of rank 29 recently discovered by Elkies and Klagsbrun. Finally we conduct some exploratory data analyses and visualizations on our dataset. We use the elliptic curve dataset from the \emph{L-functions and modular forms database} (LMFDB).
\end{abstract}


\section{Introduction}

The conjecture of Birch and Swinnerton-Dyer is a major open problem in arithmetic geometry, and is one of the six as-yet unsolved Clay Millenium prize problems \cite{carlson2006millennium}. For an elliptic curve $E$ over $\Q$, it states that the Mordell--Weil rank $r$ -- an arithmetic quantity -- is equal to the order of vanishing of the associated $L$-function $L(E,s)$ -- an analytic quantity, and furthermore gives a precise description of the leading term of the Taylor expansion of $L(E,s)$ in terms of other invariants of $E$:
\begin{equation}\label{eq:strong_bsd}
    \frac{L^{(r)}(E,1)}{r!} \overset{?}= \frac{\Omega(E/\Q) \cdot \Reg(E/\Q) \cdot |\Sha(E/\Q)| \cdot \prod_{p}c_p }{|E(\Q)_{tors}|^2}.
\end{equation}
Here, $\Omega(E/\Q)$, $\Reg(E/\Q)$ and $c_p$ denote respectively the \textsf{real period}, the \textsf{regulator}, and the \textsf{Tamagawa number at the prime $p$} of the curve; $E(\Q)_{tors}$ is the \textsf{torsion subgroup}, and $\Sha(E/\Q)$ is the \textsf{Shafarevich--Tate group} of $E$, a group currently not known (but conjectured) to be finite for all elliptic curves over $\Q$ that measures the failure of a local-global principle on the principal homogeneous spaces for $E$ (see e.g. \cite[Remark 4.1.2]{silverman2009arithmetic}). These conjectures first appear in the literature in \cite{birch1965notes}, where the authors are largely analysing the family of congruent number elliptic curves $E_D : y^2 = x^3 - Dx$; the above expression \ref{eq:strong_bsd} was given (more generally for all abelian varieties over number fields) by Tate in \cite{tate1965conjectures}. See also \cite{birch2006lieu} for a warm account by Birch of how the conjecture arose and of the many other people who contributed ideas to it. The conjecture is also remarkable as being one of the first examples (along with the Sato-Tate conjecture, now a theorem by \cite{barnet2011family}) of a relationship in pure mathematics to be discovered via experimentation and programming with an electronic computer (the EDSAC-2).

Elliptic curves naturally lend themselves to being tabulated into a database. The simplest way to do this would be via the so-called \textsf{naive height} of $E$, which is essentially\footnote{There are often constants given in the definition that we omit because we will not work with the naive height.} the height of the largest Weierstrass coefficient. It is obvious that there are finitely many elliptic curves up to bounded naive height, allowing the curves to be tabulated in order of increasing naive height. Historically however, elliptic curves have been tabulated according to the \textsf{conductor}, a more subtle invariant arising from the Galois representations associated to $E$. That there are only finitely many elliptic curves of bounded conductor follows from a theorem of Shafarevich (that there are only finitely many elliptic curves with good reduction outside of a given finite set of primes), allowing one to tabulate elliptic curves by increasing conductor. Swinnerton-Dyer was the first to do this; his table of elliptic curves up to conductor $200$ first appeared in the published literature as Table 1 in \cite{birch2006modular}; for each curve, various other features of interest are given, such as the Kodaira symbol for each prime of bad reduction, the rank, and the isogeny graph. This was expanded by Cremona \cite{cremona1997algorithms}, who extended this to conductor $1000$ via faster algorithms with modular symbols and better optimized implementations; many other curve features were also given in his tables. Subsequent editions extended this conductor bound still further, and the tables were eventually incorporated into the \emph{L-functions and modular forms database} \cite{lmfdb}, which currently has complete data for curves up to conductor $500{,}000$, as well as all curves of \emph{prime conductor} up to $300$ million.

Machine learning (ML) tools started to become ubiquitous in the mid-2010s, driven by advances in deep learning and increases in computational power (especially through GPUs). Google's releasing of TensorFlow, an open-source machine learning framework, lowered the barrier for entry into ML; around the same time, Cloud providers such as Amazon Web Services, Google Cloud, and Microsoft Azure started offering ML services, enabling companies and developers to integrate ML into applications without needing specialized hardware or expertise. By the late 2010s, with the popularity of open-source Python libraries such as scikit-learn \cite{scikit-learn} and PyTorch \cite{pytorch}, employing ML tools on datasets was commonplace. The first time this was done on a dataset of elliptic curves was in 2019, when Alessandretti, Baronchelli and He \cite{alessandretti2023machine} compared different ML models for predicting various features (including the size of $\Sha(E/\Q)$) of elliptic curves from the defining Weierstrass model. Although unsuccessful\footnote{c.f. Section 6 of \emph{loc. cit.}: ``due to the very high variation in the size of $a_i$ [the Weierstrass coefficients] ... one could not find a good machine-learning technique ... that seems to achieve this''.}, it was the first study of its kind, and paved the way for subsequent work by He, Lee and Oliver \cite{he2023curves} that also conducted several ML experiments on the elliptic curve data in the LMFDB. One of the more impressive aspects of this work was the prediction of the rank of $E$ from a limited number of trace of Frobenius values $a_p$; further analysis of this by He, Lee, Oliver and Pozdynakov \cite{he2024murmurations} led to the discovery of \emph{murmurations of elliptic curves} that remains unexplained (although see \cite{zubrilina2023murmurations} for a proof in the realm of modular forms). 

The paper \cite{he2023curves} of He--Lee--Oliver also attempted to train ML models for predicting the size of $\Sha(E/\Q)$, albeit in a more modest way than \cite{alessandretti2023machine}. Rather than attempting to train a regression model for the size itself, He--Lee--Oliver took a dataset of about $50{,}000$ curves, half of which had $|\Sha(E/\Q)| = 4$ and half had $|\Sha(E/\Q)| = 9$, and attempted to train an ML classifier on the $a_p$ values to distinguish between these cases\footnote{it is known that, if the order of $\Sha(E)$ is finite, then its order is a square, due to the alternating property of the Cassels--Tate pairing; thus, orders 4 and 9 are the smallest nontrivial values of this group.}. Despite trying a variety of methods, they were not especially successful with this task, obtaining results no better than a precision of $0.6$. Tasks of this sort are of interest because of the open conjectures about $\Sha(E/\Q)$ (chiefly its conjectured finiteness); as the discovery of murmurations attests to, successfully training ML models on mathematical data can lead to the discovery of new mathematical relationships and ideas, and in the most optimistic case for the question of predicting the order of $\Sha(E/\Q)$, could shed new light on this long-standing open conjecture.

This is the motivation for why we in the present paper continue this approach of attempting to train ML models for predicting the size of $\Sha(E/\Q)$.

\subsection{Outline of paper and summary of results}\label{ssec:main_contributions}

All of our work has used the elliptic curve database in the LMFDB. In brief, the following are the main contributions of our work.

\begin{enumerate}
    \item We improve upon the binary classification experiments of He--Lee--Oliver, obtaining $> 95\%$ accuracy in some cases. (See \Cref{sec:1-vs-4,sec:4-vs-9}.)
    \item We improve upon the regression model of Alessandretti--Baronchelli--He \cite{alessandretti2023machine}, analyze feature importances of the model, and apply the model to predict the order of $\Sha(E)$ where $E$ is the elliptic curve of rank $29$ recently discovered by Elkies and Klagsbrun \cite{ek_rank_29}. We also carry out analyses of the dataset to investigate the validity of Delaunay's heuristics on the distribution of $\Sha(E)$. (See \Cref{sec:regression_model}.)
    \item We conduct further analyses of the LMFDB dataset, motivated by attempting to understand why the models perform well, as well as potentially discovering new relationships between the BSD features. These data analyses are inconclusive; nevertheless, we feel it is still valuable to report on them. (See \Cref{sec:pca}.)
    \item We have developed and made publicly available an extensive codebase for other researchers to use and build upon. This includes several Jupyter notebook files that serve as tutorials for different parts of this paper, and explain what the code is doing. This codebase is available at:

    \begin{center}
    \url{https://github.com/barinderbanwait/ml_rnt}
    \end{center}

    Unless otherwise specified, filenames given in the paper will always be relative to the \path{experiments} directory of this repository.
    
\end{enumerate}

The binary classification experiments in \Cref{sec:1-vs-4,sec:4-vs-9} proceed by using subsets of \emph{BSD features} (the quantities that arise in \Cref{eq:strong_bsd}) for training, rather than several $a_p$ values that was originally done in He--Lee--Oliver. \Cref{sec:1-vs-4} restricts to positive rank elliptic curves and studies a binary classification problem to distinguish between $|\Sha(E)| = 1$ or $|\Sha(E)| = 4$; the motivation behind this is to force the regulator $\Reg(E)$ to be nontrivial (since for rank 0 curve the regulator is necessarily $1$). One question that motivates many of the experiments in these two sections is which of the BSD features is the most predictive of $|\Sha(E)|$; for the original He--Lee--Oliver experiment, this appears to be the real period, although for the positive rank experiment, this is the Tamagawa product. While we explore possible explanations for these observations, this remains at present mysterious and would be worthy of further study.

Finally in \Cref{sec:future} we briefly indicate possible future avenues of study.

\subsection{Acknowledgements}

This project formed at \emph{Rethinking Number Theory 5} in June 2024, an AIM Mathematical Research Community, and the authors thank the organizers of that online workshop -- Heidi Goodson, Allechar Serrano López, and Mckenzie West -- for bringing the authors of this paper together.

The project was significantly developed during the Harvard Center for Mathematical Sciences and Applications (CMSA) Program on \emph{Mathematics and Machine Learning} which took place in September and October 2024, and particularly during the two number theory weeks that allowed three of us (AB, BSB, XH) to meet and work in person. We thank Edgar Costa, Michael Douglas, and Andrew Sutherland for organising these activities and for putting in place computational resources including GPUs that were instrumental in our work.

We are also grateful to the organizers of the workshop \emph{Murmurations in Arithmetic Geometry and Related Topics} -- Yang-Hui He, Abhiram Kidambi, Kyu-Hwan Lee, and Thomas Oliver -- held at the Simons Center for Geometry and Physics in Stony Brook, NY, that again allowed for three of us (AB, BSB, XH) to meet, develop, and present this work.

BSB acknowledges support from the Simons Foundation, grant \#550023 for the Collaboration on Arithmetic Geometry, Number Theory, and Computation.
AF is supported by a Croucher Scholarship.

\section{Binary classification between orders 4 and 9}
\label{sec:4-vs-9}

Since there is already a set benchmark in the literature from the work of He--Lee--Oliver, we take their classification problem as our point of departure and attempt to obtain a model with higher accuracy. We therefore in this section limit ourselves to working with essentially the same dataset as they did in Section 4.5 of their paper \cite{he2023curves}. Specifically,  we use a balanced subset of curves in the LMFDB of conductor up to $500,000$ and $|\Sha (E/\Q)|$ equal to $4$ and $9$. 
\begin{table}[h!]
\centering
\begin{tabular}{c|c}

\(|\Sha(E/\mathbb{Q})|\) & \# curves \\
\hline
4 & 50428 \\
9 & 50428 \\

\end{tabular}
\end{table}

The dataset is then shuffled, and we reserve a random subset of $20\%$ of curves as a test set, kept unseen during training. We fix the random seed for reproducibility.

One must necessarily use more features than merely the $a_p$ values, since these are isogeny invariant, while $\Sha(E/\Q)$ is not\footnote{An easy search in the LMFDB reveals plenty of isogenous curves with differing order of $\Sha$.}. Rearranging \Cref{eq:strong_bsd} one obtains

\begin{equation}\label{eq:bsd_for_sha}
    |\Sha(E/\Q)| \overset{?}= \frac{|E(\Q)_{tors}|^2 \cdot L^{(r)}(E,1)/r!}{\Omega(E/\Q) \cdot \Reg(E/\Q) \cdot \prod_{p}c_p}.
\end{equation}

Therefore, in our experiments, we primarily consider what we call the \textit{BSD features}: the special value $L^{(r)}(E,1)/r!$, the torsion size $|E(\Q)_{tors}|$, the real period $\Omega(E/\Q)$, the regulator $\Reg(E/\Q)$, and the Tamagawa product $\prod_{p}c_p$.

\subsection{Using all BSD features} \label{subsec:allbsd}

Assuming equation~\ref{eq:bsd_for_sha} to be true, one would expect a model trained with \emph{all} of the features on the right-hand side of this equation to perform rather well, particularly if it is a model designed to detect multiplicative relationships among the features. This is our first experiment, and is confirmed to be the case in \path{balanced_4_9_all_BSD.ipynb}. In this file:

\begin{enumerate}
    \item We train a logistic regression classification model on all features, and obtain 64\% accuracy. This is provided merely as a benchmark.

    \item We train a logistic regression classification model on the \emph{logarithm} of all features, and obtain perfect 100\% accuracy.

    Since logistic regression is designed to find linear relations among features, and since taking logarithms of both sides of \Cref{eq:bsd_for_sha} yields a linear relation between the features, it is expected that logistic regression should perform very well.

    \item We train an ordinary least squares linear regression model to detect the linear relation itself, confirming the validity of \Cref{eq:bsd_for_sha}.

    We do not, however, claim that this is in any way evidence for BSD, since the LMFDB \path{sha_order} field is computed via \Cref{eq:bsd_for_sha}.

    \item We train a histogram-based gradient-boosting classification tree on all features and obtain 98\% accuracy.

    We use this classifier due to its ease of use. Models using this classifier are more flexible than logistic regression in detecting relationships between features, including products and ratios, since they sequentially improve predictions by learning residuals. It makes no difference whether or not we take the logarithm of the data (as is verified in the notebook).
\end{enumerate}

\subsection{Removing one BSD feature at a time}
\label{subs: remove}

Clearly, training a model with all terms in the conjectural BSD formula is expected to yield high accuracies. It would be more valuable to investigate whether one can train a model that can successfully predict the size of $\Sha$ with \emph{incomplete} information, particularly information that is easy to obtain. In particular, we ask which of the above five BSD features are most predictive of the size of $\Sha$.

\begin{figure}[htbp] 
    \centering
    \includegraphics[scale=0.5]{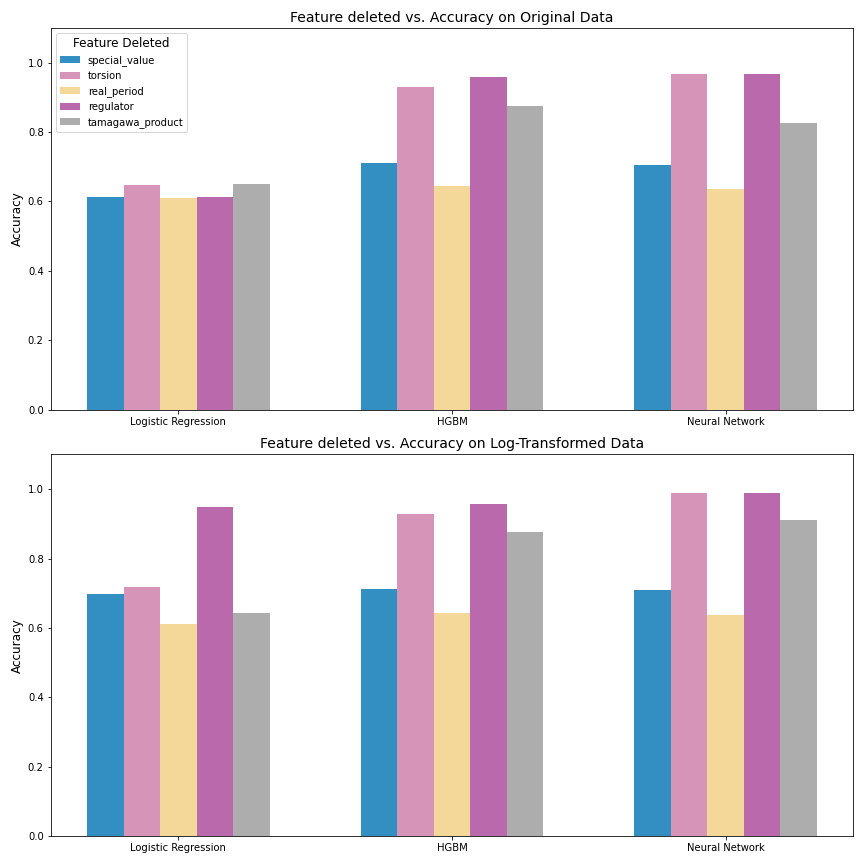}
    \caption{Feature Deleted vs Accuracy Across Models for $|\Sha(E/\Q)| = 4$ and $| \Sha(E/\Q)| = 9$.}
    \label{fig: feat_acc}
\end{figure}

\begin{experiment}
\label{exp: 4-vs-9}

We compare the performance of three models on the data. We use the logistic regression and the tree models described above. The success of the tree model in the previous section motivates us also to try a neural network with one hidden layer.


Specifically, we consider the following binary classification models, using the original and log-transformed data:

\begin{enumerate}
\item logistic regression (in \path{balanced_4_9_logistic_remove_one.ipynb});
\item a histogram-based gradient boosting classification tree (in \path{balanced_4_9_hgbm_remove_one.ipynb}); and
\item a feedforward neural network with 3 hidden layers with 128, 64, and 32 hidden units  (in \path{balanced_4_9_NN_remove_one.ipynb}). The model uses ReLU activation functions, dropout of $0.3$, cross-entropy loss function, Adam optimizer and a learning rate of $0.001$. We run the model for 100 epochs and record the best test accuracy. 
\end{enumerate}

    For each model, we remove in turn one of the five features -- the special value, the torsion, the real period, the regulator and the Tamagawa product -- and record the accuracy on the test set.
    The lower the resulting accuracy scores, the more important we expect the feature to be. The performance of the models is illustrated in Figure \ref{fig: feat_acc}.
\end{experiment}


The performance of the logistic regression classification model on the original data is relatively poor, and removing individual features seems to make little difference to the outcome. This aligns with the observation made in \Cref{subsec:allbsd} that logistic regression finds linear relations among features. Indeed, even after log-transforming the data, removing any feature results in a significant drop in performance with the exception of the regulator, which is briefly discussed below.

We also notice that log-transforming the data seems to generally improve performance. This gives the most obvious improvement in the logistic regression model (when removing the regulator) and slight improvements in the neural network. However, the performance of the tree model is not affected by log-transforming the data at all, which is consistent with our expectation that transforming the data monotonically should not affect the performance.

For the high-performing models, which include the tree models and neural networks, the real period and the special value seem to be the most important features, and neither model is able to make up for the loss of those. They are followed by the Tamagawa product, and finally the torsion and regulator. We suspect that the reason why the regulator is among the least important features is that a vast majority ($92.6\%$) of the curves in the dataset have rank $0$, and therefore trivial regulator, which makes this feature nearly constant. We perform a similar experiment on positive rank curves in \Cref{exp: 1-vs-4}. 

\subsection{Training with $a_p$ values}

\begin{figure}[htbp] 
    \centering
    \includegraphics[width=\linewidth]{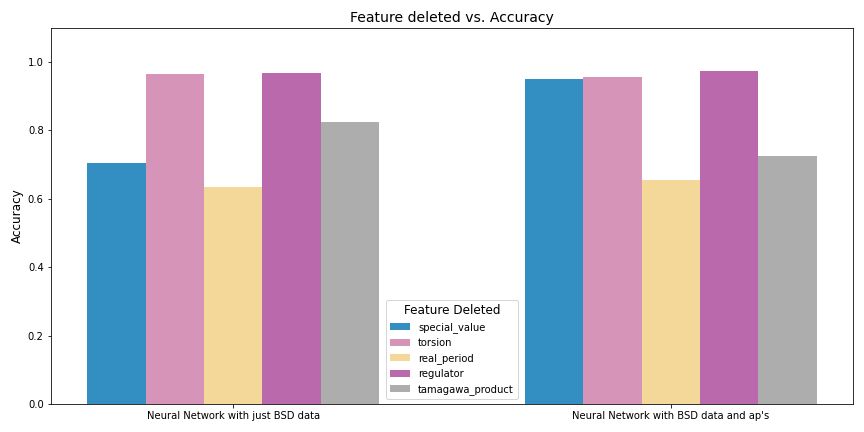}
    \caption{Feature deleted vs accuracy in a feedforward neural network classification task between between: without and with $a_p$ values. }
    \label{fig: aps_class}
\end{figure}

The features in the original experiments of He, Lee and Oliver \cite{he2023curves} consisted of $500$ $a_p$ values. While we don't expect these values to contribute to training in addition to the BSD features in our earlier experiments, we verify this using the neural network experiment done in \Cref{subs: remove}. Namely, we compare the performance of the classifier when we use the BSD features without log-transforming to that when we also add the first $100$ $a_p$ values in \path{balanced_4_9_aps.ipynb}. In both cases, we remove one BSD feature at a time and record the results in Figure \ref{fig: aps_class}.

In general, adding the $a_p$ values does not improve the model significantly. The only exception is when removing the special value; this is to be expected since the $a_p$ values determine the $L$-function of the curve and thus already encode the special value.

\section{Binary classification on positive rank curves between orders 1 and 4}
\label{sec:1-vs-4}

In \Cref{sec:4-vs-9}, we performed experiments on the dataset in \cite{he2023curves}. As mentioned there, that dataset contains many curves of rank $0$, which therefore have trivial regulator. In this section, we perform similar experiments on a dataset of positive rank curves.  

Specifically, we use a balanced set of curves from the LMFDB of conductor up to $500,000$, rank $\rk(E)>0$ and  $|\Sha (E/\Q)|$ equal to $1$ and $4$. We restrict ourselves to $|\Sha(E/\Q)| \in \{1, 4\}$ because at the time of this article, the LMFDB dataset contains only 1462 positive rank curves with $|\Sha (E/\Q)| = 9$. 

\begin{table}[h!]
\centering
\begin{tabular}{c|c}

\(|\Sha(E/\mathbb{Q})|\) & \# curves \\
\hline
1 & 18710 \\
4 & 18710 \\

\end{tabular}
\end{table}

The dataset is shuffled and we reserve a random subset of $20\%$ of curves as a test set, kept unseen during training. We fix the random seed for reproducibility.

\begin{figure}[htbp] 
    \centering
    \includegraphics[width=\linewidth]{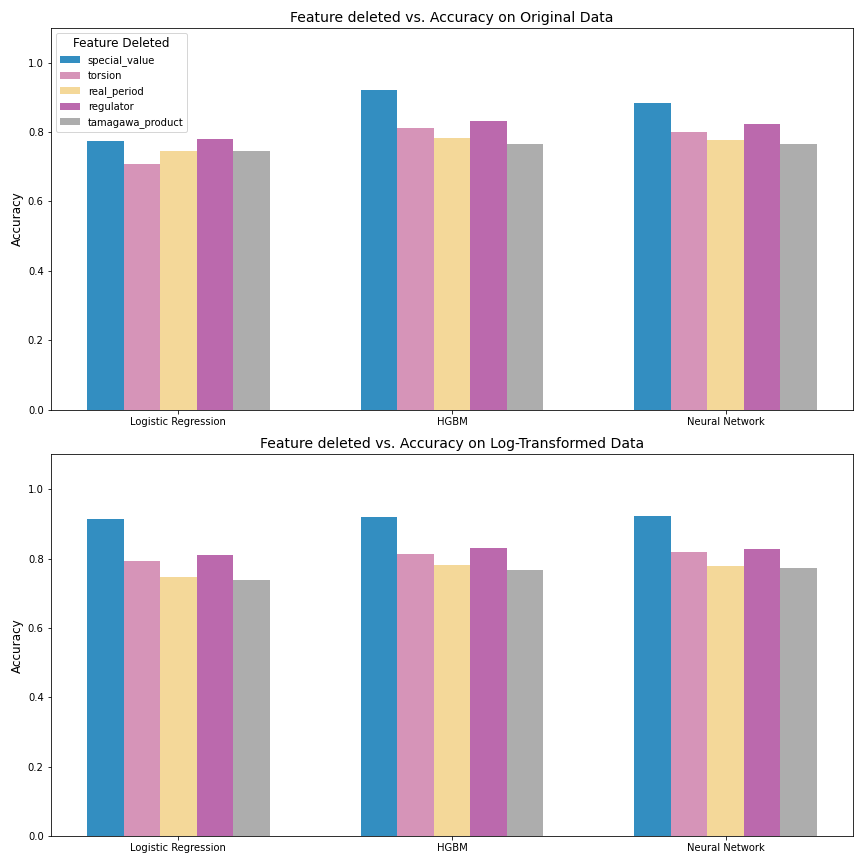}
    \caption{Feature Deleted vs Accuracy Across Models for $|\Sha(E/\Q)| = 1$ and $|\Sha(E/\Q)| = 4$.}
    \label{fig: feat_acc14}
\end{figure}

\begin{experiment}
\label{exp: 1-vs-4}

We compare the performance of the three binary classification models (logistic regression in \path{balanced_1_4_logistic_remove_one.ipynb}, a tree model in \path{balanced_1_4_hgbm_remove_one.ipynb}, and a feedforward neural network in \path{balanced_1_4_NN_remove_one.ipynb}) on the original data and the log-transformed data. We use the same setup as in \Cref{exp: 4-vs-9}. For each model, we remove in turn one of the five features  $L^{(r)}(E,1)/r!$, $|E(\Q)_{tors}|$, $\Omega(E/\Q)$, $\Reg(E/\Q)$,  and $\prod_{p}c_p$, and record the accuracy on the test set. We record their performance in Figure \ref{fig: feat_acc14}. 
\end{experiment}

We notice here that log-transforming the data improves the performance of the logistic regression model significantly. Moreover, upon log-transforming the data, logistic regression performs similarly to the non-linear models.

The real period, as in \Cref{exp: 4-vs-9}, is again among the most predictive features. Curiously, the Tamagawa product also seems to be highly indicative of $|\Sha(E/\Q)|$, and the special value is now the least predictive feature. We explore the relationships between these features in more detail in \Cref{sec:pca}.

\section{Regression for Sha size}\label{sec:regression_model}
The models in this paper so far have been focused on binary classification. In this section we explain our steps for training a regression model to predict the size of $\Sha$ from various other BSD features. 

\subsection{Feature And Model Selections}
We trained our models on an $80\%$ sample of the curves in the LMFDB of conductor up to $500{,}000$, and tested them both on the remaining $20\%$ of this dataset, \emph{as well as} the curves of prime conductor between $500{,}000$ and $300$ million. This was to investigate whether a model trained on curves of bounded conductor could reasonably predict the order of $\Sha$ of curves beyond what the model had been trained on. Since one application we have in mind of the usefulness of such models is their ability to predict the order of $\Sha$ for curves of large rank (where the conductor is necessarily large), we were curious to know how the model would fare on this larger conductor dataset.

The dataset we obtained from the LMFDB had, in addition to the BSD features in \Cref{eq:bsd_for_sha}, also the first $100$ $a_p$ values, and the following features: rank, conductor, adelic level, adelic index, adelic genus, and the PARI-encoded Kodaira symbols.

We attempted a variety of different models, and experimented with different neural network architectures. A key challenge was that approximately \(90\%\) of the dataset consists of curves with trivial \( \Sha \), so finding a model that does not overfit the data is particularly challenging. Furthermore, a naive approach that predicts all curves have trivial \( \Sha \) already achieves about \(90\%\) accuracy, leaving little room for improvement, at least in terms of accuracy. In particular, our neural network models did not show significantly better performance than the naive approach. In our experiments, the histogram-based gradient boosting machine demonstrated better performance compared to the other models discussed in Section \ref{exp: 1-vs-4}, leading us to select it as our model of choice.

\begin{figure}[h]
    \centering
    \includegraphics[width=0.8\linewidth]{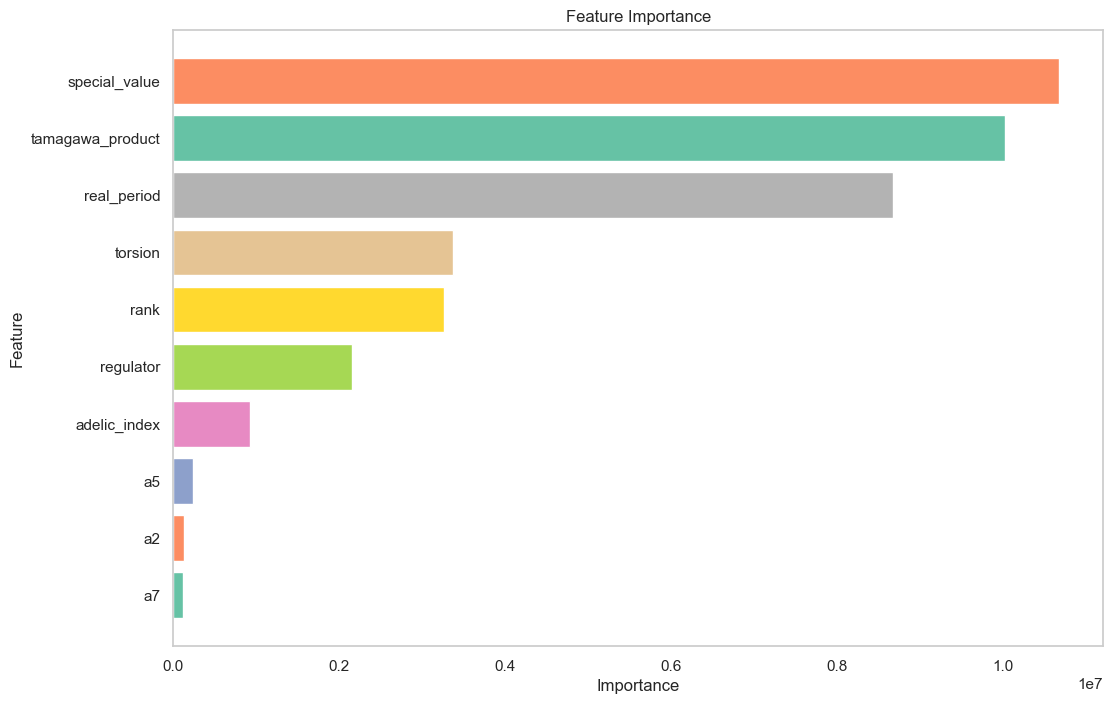}
    \caption{The importance of the 10 most significant features computed using LightGBM. The values represent the information gain contributed by each feature to the model.}
    \label{fig: lightgbm_feature}
\end{figure}

As for feature selections, we started by training a regression model using LightGBM \cite{ke2017lightgbm}, a variant of the histogram-based gradient boosting machine, on all of the features available. This was done to get a sense of what features are important for $|\Sha(E/\Q)|$ predictions, since this model can quantify feature importances. The importance of the 10 most significant features is shown in \Cref{fig: lightgbm_feature}. 

Given the importance of rank described in \Cref{fig: lightgbm_feature}, we subsequently conducted three experiments using a histogram-based gradient boosting machine as the model: 
\begin{itemize}
    \item substituting \(\Reg(E/\Q)\) with \(r\);
    \item a baseline experiment excluding both \(\Reg(E/\Q)\) and \(r\);
    \item an experiment using all BSD features as a benchmark for comparison.
\end{itemize}

These experiments are carried out in \path{tree_regression_model.ipynb}. 

\subsection{Results}{\label{subsec:reg_results}}
Since \( \sqrt{|\Sha(E/\mathbb{Q})|} \in \mathbb{Z}^+ \), we train the model to predict \( \sqrt{|\Sha(E/\Q)|} \), then rounding the regression predictions to the nearest positive integer. The model's performance is then evaluated using accuracy score as well as the Matthews Correlation Coefficient (MCC)
\footnote{
The Matthews Correlation Coefficient (MCC) is defined as:

\begin{equation}
\text{MCC} = \frac{\text{TP} \cdot \text{TN} - \text{FP} \cdot \text{FN}}{\sqrt{(\text{TP} + \text{FP})(\text{TP} + \text{FN})(\text{TN} + \text{FP})(\text{TN} + \text{FN})}},
\end{equation}
which ranges from $-1$ to $+1$, where $+1$ indicates perfect predictions, and $-1$ indicates predictions are completely opposite to the true values.
}. 
We include MCC as it is a more robust metric for imbalanced datasets. In particular, the naive approach of predicting that all curves have trivial \( \Sha \) has an MCC of \( 0 \).

The accuracy scores and MCC values for all experiments are presented in \Cref{tab:regression_model_res}. These results demonstrate that the model achieves high accuracy and MCC overall, while also generalizing effectively to datasets of curves with larger conductors, though with a minor decrease in performance compared to the small conductor dataset.

\begin{table}[h!]
\centering
\begin{tabular}{lcccc}
\hline
\textbf{Features}                  & \multicolumn{2}{c}{\textbf{Small Conductor}} & \multicolumn{2}{c}{\textbf{Large Conductor}} \\ 
                                   & \textbf{Accuracy} & \textbf{MCC}            & \textbf{Accuracy} & \textbf{MCC}            \\ \hline
BSD features                       & 0.99             & 0.95                    & 0.97             & 0.86                    \\ 
Substituting \(\Reg(E/\Q)\) with \(r\) & 0.99             & 0.94                    & 0.98             & 0.87                    \\ 
Excluding \(\Reg(E/\Q)\) and \(r\) & 0.91             & 0.38                    & 0.89             & 0.22                    \\ 
\hline
\end{tabular}
\caption{Comparison of Accuracy and MCC for Small and Large conductor sets. `Small conductor' here means conductor at most $500{,}000$; `Large conductor' means prime conductor between $500{,}000$ and $300$ million.}
\label{tab:regression_model_res}
\end{table}

When all BSD features are included, the model achieves an accuracy of 0.97 and an MCC of 0.86 on the large conductor dataset, compared to 0.99 accuracy and 0.95 MCC for the small conductor dataset. Similarly, when substituting \(\Reg(E/\Q)\) with \(r\), the performance remains robust, with an accuracy of 0.98 and an MCC of 0.87 for large conductors, which are slightly higher than those achieved using all BSD features. However, excluding both \(\Reg(E/\Q)\) and \(r\) results in a significant drop in performance, with the accuracy and MCC falling to 0.89 and 0.22, respectively, for large conductors. 

\begin{figure}[h]
\centering
\begin{subfigure}{0.49\linewidth}
    \centering
    \includegraphics[width=\linewidth]{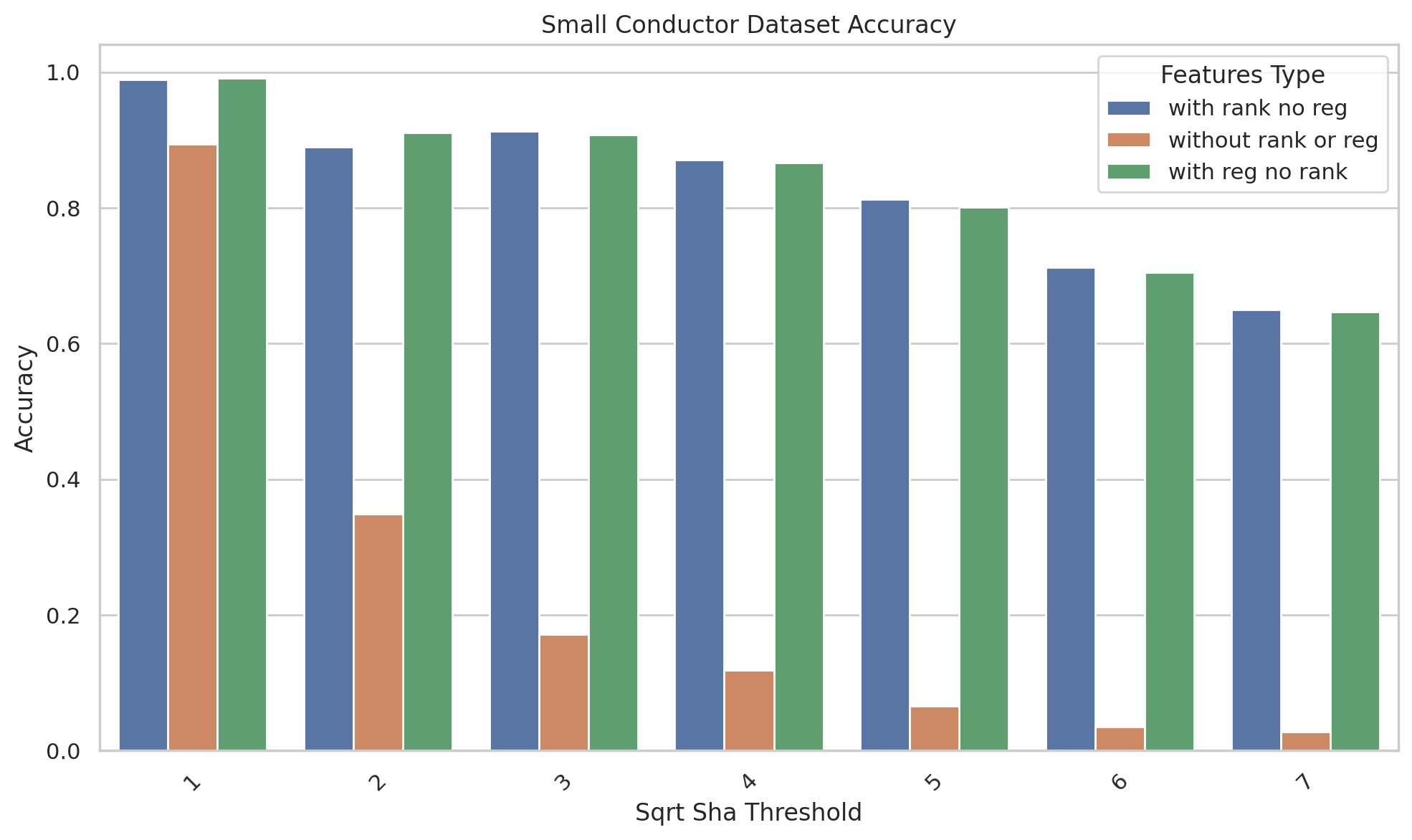}
    \caption{The small conductor case.}
\end{subfigure}
\begin{subfigure}{0.49\linewidth}
    \centering
    \includegraphics[width=\linewidth]{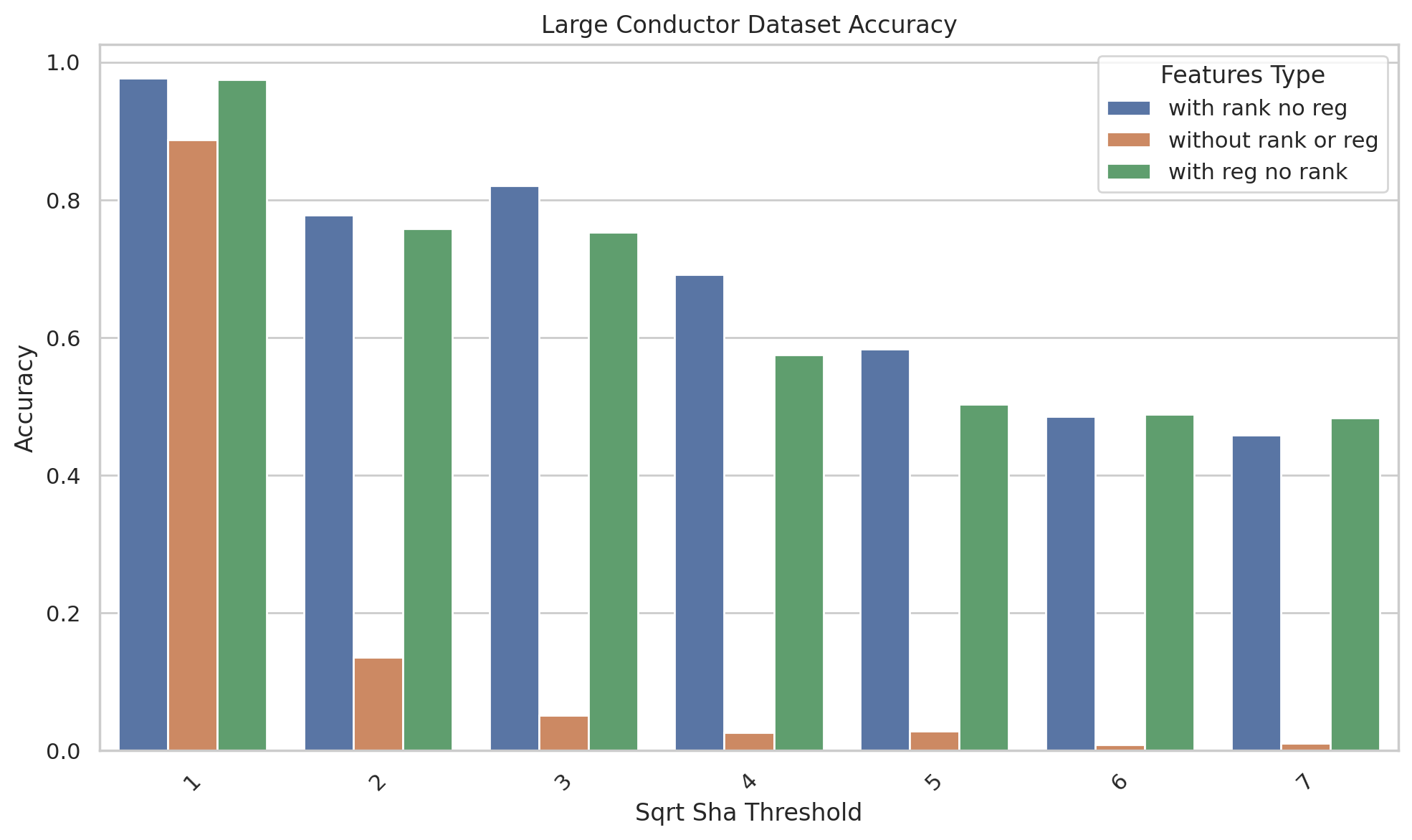}
    \caption{The big conductor case.}
\end{subfigure}
\caption{The accuracy within subsets of curves where \(\sqrt{\Sha} \geq \) a given threshold for both the small conductor and large conductor datasets. The results are comparable between the model that includes all variables in the BSD formula and the model that substitutes \(\Reg(E/\Q)\) with \(r\). The model that excludes both variables performs significantly worse.}
\label{fig:rank_reg_same}
\end{figure}

The results above suggest that replacing the regulator with the rank during training has no significant impact on performance. To explore this further, we analyze the prediction accuracy of the models when restricted to subsets of curves with a minimum threshold on \( |\Sha(E/\Q)| \). The results are presented in \Cref{fig:rank_reg_same}.

However, the conclusion that replacing the regulator with the rank has no significant impact on performances when predicting $|\Sha(E/\Q)|$ cannot be drawn. In particular, 
when the models are separately trained on $r > 0$ curves, $\Reg(E/\Q)$ has significantly more prediction power than $r$. When the models are separately trained on $r = 0$ curves, they have the same performances mostly likely due to the model being able to infer that the $r = 0$ implies that $\Reg(E/\Q) = 1$.

\begin{figure}[h]
    \centering
    \begin{subfigure}[t]{0.49\linewidth}
        \centering
        \includegraphics[width=\linewidth]{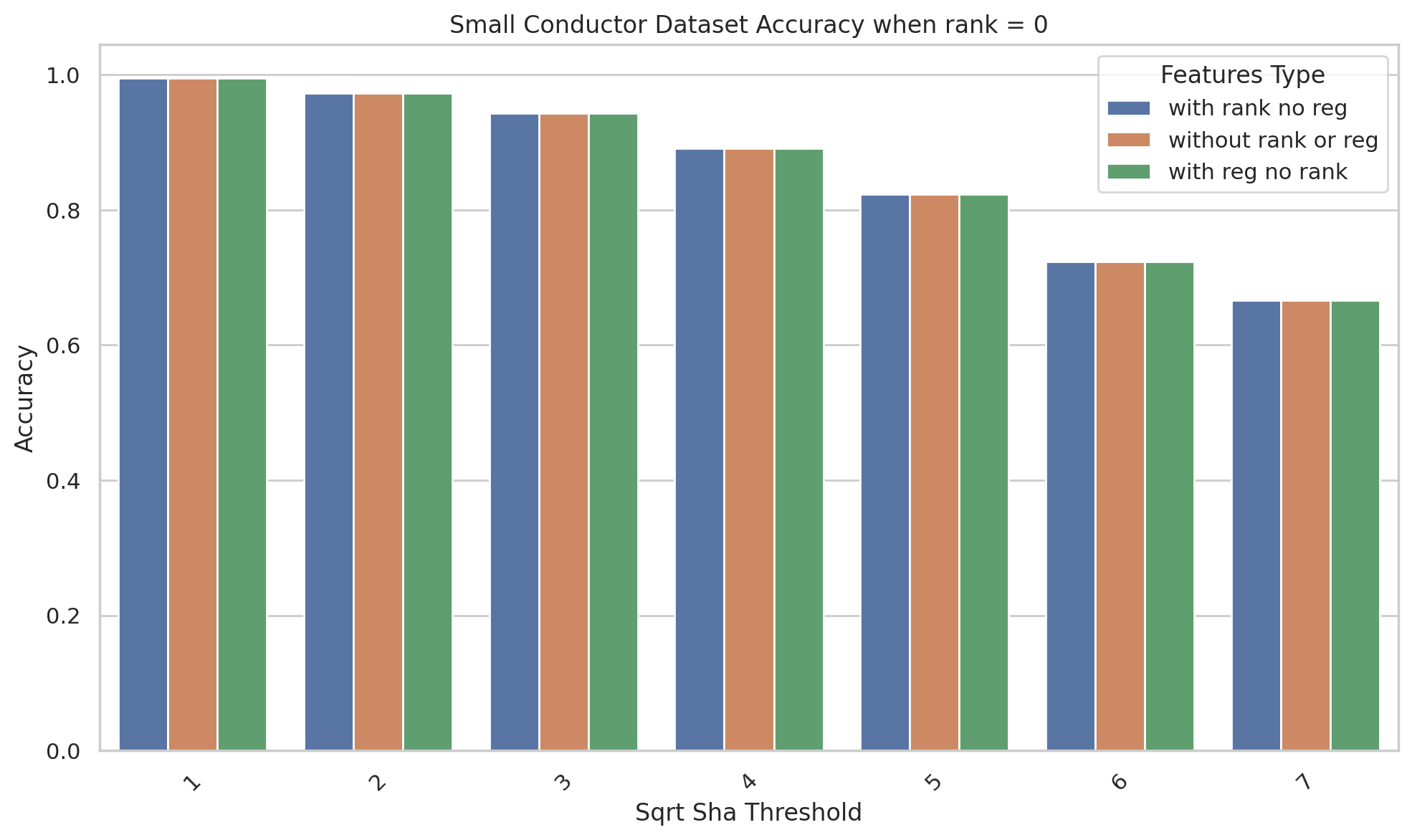}
        \caption{Accuracy when $r = 0$}
        \label{fig:acc_rank_0}
    \end{subfigure}
    \hfill
    \begin{subfigure}[t]{0.49\linewidth}
        \centering
        \includegraphics[width=\linewidth]{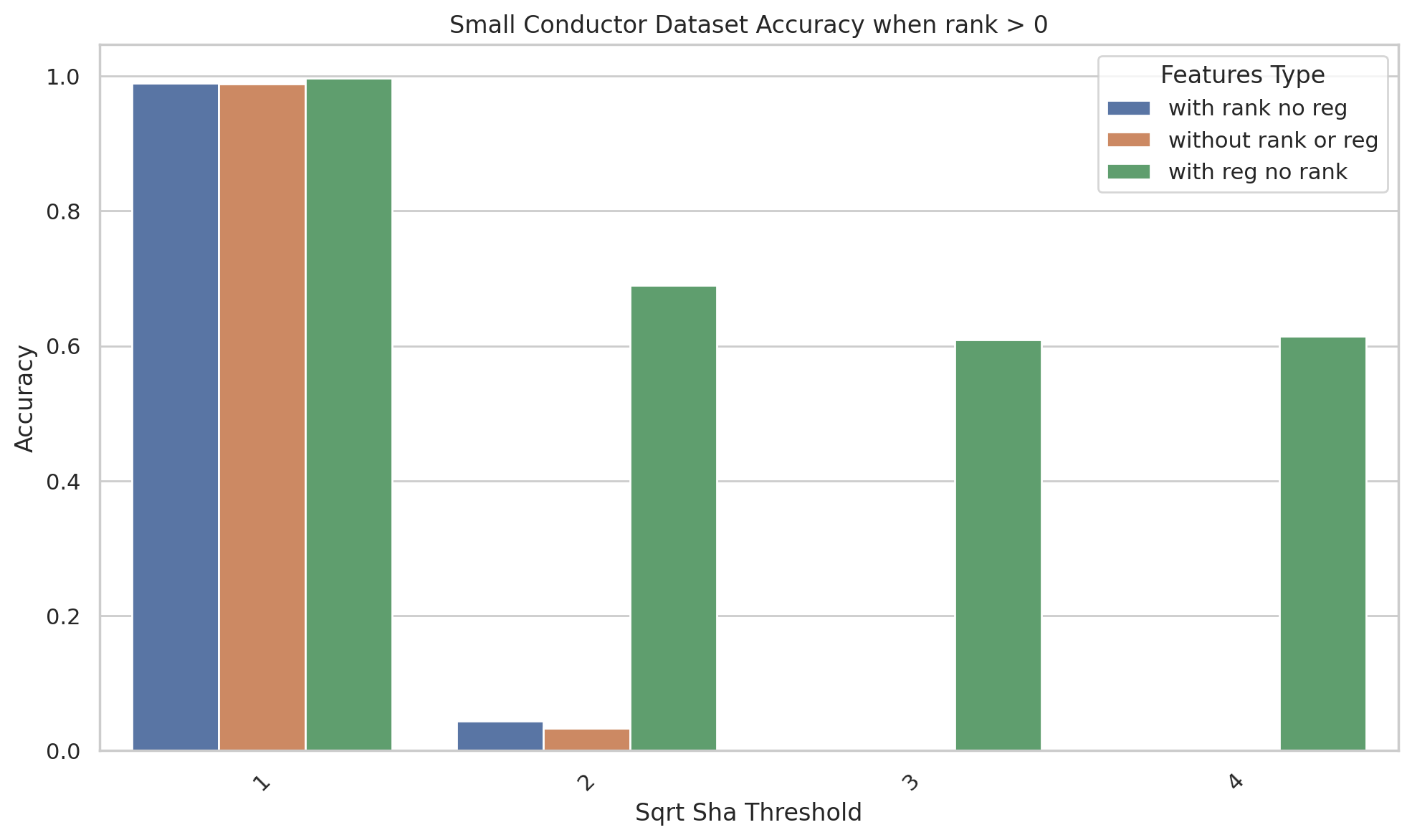}
        \caption{Accuracy when $r > 0$}
        \label{fig:acc_rank_gt_0}
    \end{subfigure}
    \caption{Comparison of accuracy between models when training with $r = 0$ and $r > 0$ curves separately. While no differences are observed between the models for rank $0$ curves, including \(\Reg(E/\Q)\) is essential for accurately predicting $|\Sha(E/\Q)|$.}
    \label{fig:accuracy_comparison}
\end{figure}

Finally, we emphasize that the results presented above were achieved without oversampling the majority class (curves with \( |\Sha(E/\Q)| = 1 \)) or undersampling the minority class (curves with \( |\Sha(E/\Q)| > 1 \)), despite the inherent imbalance in the dataset. This approach was intentional, as our aim is to develop a methodology that can accurately predict \( |\Sha(E/\Q)|\), including cases such as an extreme example curve in which certain values from \Cref{eq:bsd_for_sha} are challenging to compute. For a concrete example, see \Cref{subsec:Elkies}, where we apply a trained model on the elliptic curve with rank 29 recently discovered by Elkies and Klagsbrun \cite{ek_rank_29}.

\subsection{On Delaunay's Heuristics}
In analogy to the Cohen-Lenstra heuristics for class groups \cite{CL84}, Delaunay made a precise conjecture on the distribution on the structure of $|\Sha(E/\Q)|$ in \cite{delaunay2001heuristics, delaunay2007heuristics, delaunay2013pelltorsionpointsfiniteabelian} when $r = 0$ and $r = 1$, which is later proved to be compatible with the conjectured model in \cite{bhargava2013modeling}. As a consequence of the conjectures, $\Sha(E/\Q)$ is expected to be ``small'' when $r > 0$, and ``large" when $r = 0$. Specifically, we expect that when $r =0$, the probability that $\Sha(E/\Q)$ is isomorphic to the square of a cyclic group is approximately 0.977076, and when $r = 1$, the probability that $|\Sha(E/\Q)| = 1$ is approximately 0.54914 according to \cite{delaunay2001heuristics}. 

\begin{figure}[h]{\label{fig:prop of sha sizes}}
\centering
\begin{subfigure}{.5\textwidth}
  \centering
  \includegraphics[width=\linewidth]{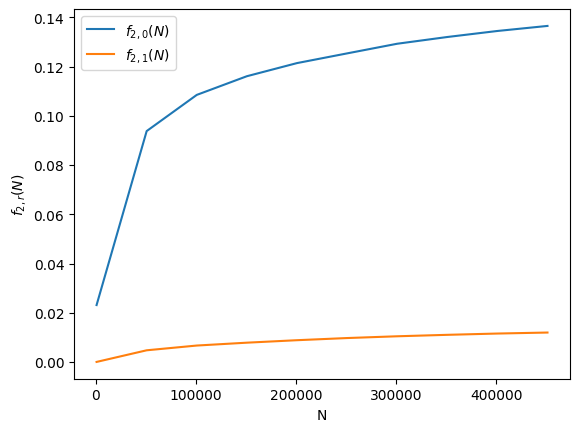}
  \caption{$p = 2$}
  \label{fig:sub1}
\end{subfigure}%
\begin{subfigure}{.5\textwidth}
  \centering
  \includegraphics[width=\linewidth]{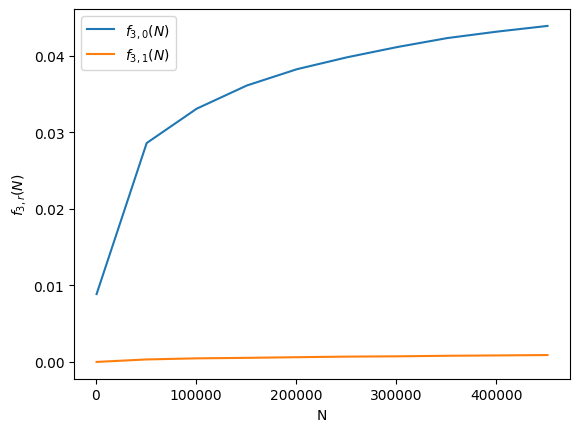}
  \caption{$p = 3$}
  \label{fig:sub2}
\end{subfigure}
\caption{The proportion of $p \mid | \Sha(E/\Q)|$  curves up to conductor $N$ when $r = 0$ and $r = 1$. $f_{r,p}(N)$ is the proportion of curves with $|\Sha(E/\Q)|$ divisible by $p$ within rank $r$ ones up to conductor $N$.}
\label{fig:p_div_sha}
\end{figure}

\begin{figure}[h]
    \centering
    \includegraphics[width=0.5\linewidth]{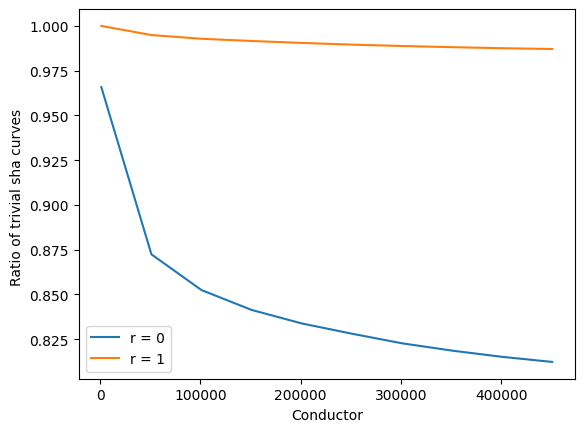}
    \caption{The proportion of $|\Sha(E/\Q)| = 1$  curves up to conductor $N$ when $r = 0$ and $r = 1$. }
    \label{fig: trivial_sha}
\end{figure}

Our experiment provides partial evidence of Delaunay's Heuristic by the model's recognition of the contribution to $|\Sha(E/\Q)|$. Moreover, we conduct further analysis of the empirical distributions of $|\Sha(E/\Q)|$ using the dataset from LMFDB, which includes all curves with conductors up to 500,000, and compute 
$$
f_{p, r}(N) = \mu \bigg(p \mid |\Sha(E/\Q)| \bigg|  \text{conductor}(E) < N, \text{rank}(E) = r \bigg),
$$
where $\mu$ here is the empirical measure. For $p = 2$ and $3$, the results can be found in \Cref{fig:p_div_sha}. We also compute the proportion of \( |\Sha(E/\Q)| = 1 \) curves up to conductor \( N \) in \Cref{fig: trivial_sha}. 

These proportions computed from the data are not close to the ones conjectured in \cite{delaunay2001heuristics}; see \Cref{tb:Delaunay-vs-observed}. However, due to the trends of the curves in \Cref{fig:prop of sha sizes} and \Cref{fig:p_div_sha}, we could expect them to become closer as the maximum conductor increases.

\begin{table}[h]
\centering
\begin{tabularx}{\textwidth}{lXXXX}
\hline
                   & \multicolumn{2}{c}{\textbf{Rank 0}}      & \multicolumn{2}{c}{\textbf{Rank 1}}       \\ \hline
                   & \textbf{Delaunay's Heuristic} & \textbf{Observed}   & \textbf{Delaunay's Heuristic} & \textbf{Observed}    \\ \hline
$|\Sha(E/\Q)| = 1$       & $0.022924$         & $0.809611$ & $0.54914$          & $0.986610$  \\
2 divides $|\Sha(E/\Q)|$ & $0.580577$         & $0.138529$ & $0.31146$          & $0.012370$  \\
3 divides $|\Sha(E/\Q)|$ & $0.360995$         & $0.044557$ & $0.0416$           & $0.0004953$ \\ \hline
\end{tabularx}
\caption{The comparison of proportions of curves with different divisibility properties among rank \(0\) and rank \(1\) curves. The values are based on the ones provided in \cite{delaunay2001heuristics} for Delaunay's heuristic, while the observed values are computed using all curves with conductors less than 500,000.}
\label{tb:Delaunay-vs-observed}
\end{table}

\subsection{Application to the rank 29 Elkies--Klagsbrun curve}\label{subsec:Elkies}
The motivation for developing a regression model is to predict $|\Sha(E/\Q)|$ for an unknown curve where some quantities in \Cref{eq:bsd_for_sha} are computationally expensive. Unlike classification models, a trained regression model can theoretically predict any value for $|\Sha(E/\Q)|$ of a given curve. As a proof of concept for the potential applications of such models, we use it to predict $|\Sha(E/\Q)|$ for the recently discovered record-breaking curve $E_{29}$ with rank 29 (under GRH) by Elkies and Klagsbrun \cite{ek_rank_29}. 

The torsion of $E_{29}$ is already known to be trivial. In his announcement, Elkies gives the gp code to compute the regulator of $E_{29}$. We subsequently computed the real period, and Tamagawa product of $E_{29}$ and summarize the information of $E_{29}$ in \Cref{tb:Elkies_BSD_values}.

\begin{table}[h!]
\centering
\begin{tabular}{cc}
\hline
\multicolumn{1}{c}{\textbf{Variable}} & \multicolumn{1}{c}{\textbf{Value}} \\ \hline
Rank                & \(29\) \\
Regulator           & \(433744182671713097629179252379019849.493842\) \\
Torsion             & \(1\) \\
Real Period         & \(3.5090427060633614999186666781786131525 \times 10^{-15}\) \\
Tamagawa Product    & \(10725120\) \\ 
\hline
\end{tabular}
\caption{Values of the variables in the BSD formula and the rank for $E_{29}$, except for the special value.}
\label{tb:Elkies_BSD_values}
\end{table}

However, computing the special value requires computing the $a_n$ trace of Frobenius coefficients for approximately $n \leq \sqrt{\mbox{conductor}}$, which in this case is approximately $4 \times 10^{74}$ terms, making it infeasible for all current computer algebra systems to calculate.

We therefore trained a histogram-based gradient boosting machine model on the rank, together with the BSD features \emph{except} the special value, since in this case we cannot compute it. As before we used 80\% of randomly selected curves from all currently available elliptic curves over $\Q$ in the LMFDB, while reserving the remaining 20\% as a test set to evaluate the model. This is carried out in \path{regression_model_Elkies_Klagsbrun.ipynb}. The dataset that contains all currently available elliptic curves over $\Q$ includes all curves with conductors up to 500,000, as well as those with prime conductors up to 300 million. 

The model predicts $E_{29}$ to have a trivial $\Sha$. It is important to note that the special value is a critical feature for determining $|\Sha(E/\Q)|$ (cf. \Cref{sec:4-vs-9}), which limits the model's performance compared to those discussed in \cref{subsec:reg_results}. On the test set, this model achieves an accuracy of 0.905 and a Matthews correlation coefficient (MCC) of 0.360.

This is unsurprising if we assume that most curves of positive rank have a trivial $\Sha$, which aligns with the patterns observed in our dataset. Additionally, to explore how machine learning models assess the likelihood of $E$ having a trivial $\Sha$, we trained separate classification models to predict this probability using the same set of features. Both the neural network model and the histogram-based gradient boosting machine predicted a probability of 
$1$ for $E_{29}$ having a trivial $\Sha$. The implementations of these models are available in \path{classification_model_Elkies_Klagsbrun.ipynb}.

\section{PCA and visualizations}\label{sec:pca}

Having conducted various ML experiments with the LMFDB dataset relating to the size of $\Sha(E/\Q)$, in this section we carry out some further analyses of the data that are prompted by the previous experiments. The code for this section is available in \path{pca_entire_dataset.ipynb}.

\subsection{PCA for the He--Lee--Oliver dataset}
We conduct a Principal Component Analysis (PCA) on the dataset used for \Cref{exp: 4-vs-9}.

The dataset is log-transformed and thereafter
normalized for each feature to have mean $0$ and standard deviation $1$. As in \Cref{sec:4-vs-9}, this dataset has 50,428 curves each of $|\Sha(E/\Q)| = 4$ and $|\Sha(E/\Q)| = 9$.

\begin{figure}[htbp] 
    \centering
    \includegraphics[width=\linewidth]{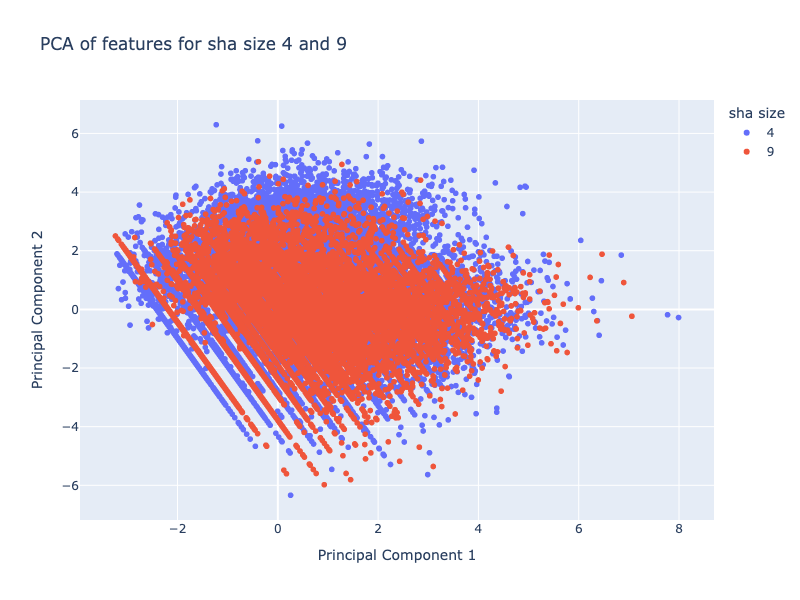}
    \caption{PCA on the dataset used in \Cref{sec:4-vs-9}.}
    \label{fig:pca_sha_4_vs_9}
\end{figure}

The first two principal components carry about 36\% and 28\% of the variance respectively, but these components are not effective at separating the two classes of curves as depicted in figure \Cref{fig:pca_sha_4_vs_9}.

Loadings for the first two principal
components are given in \Cref{tb:pc_loadings_4_9}.

\begin{table}[h!]
\centering
\begin{tabular}{lcc}
\hline
\textbf{Feature}         & \textbf{PC1} & \textbf{PC2} \\
\hline
tamagawa\_product & 0.70 &  0.12 \\
special\_value    &  -0.09 &  0.72 \\
torsion           & 0.51 &  0.29 \\
regulator         & 0.06 &  0.48 \\
real\_period      &  0.49 & 0.39 \\
\hline
\end{tabular}
\caption{PCA loadings for the first two principal components with $\Sha$ size 4 vs 9.}
\label{tb:pc_loadings_4_9}
\end{table}

\subsection{PCA for positive rank, $\Sha$ size $1$ or $4$ dataset}

\begin{figure}[htbp] 
    \centering
    \includegraphics[width=\linewidth]{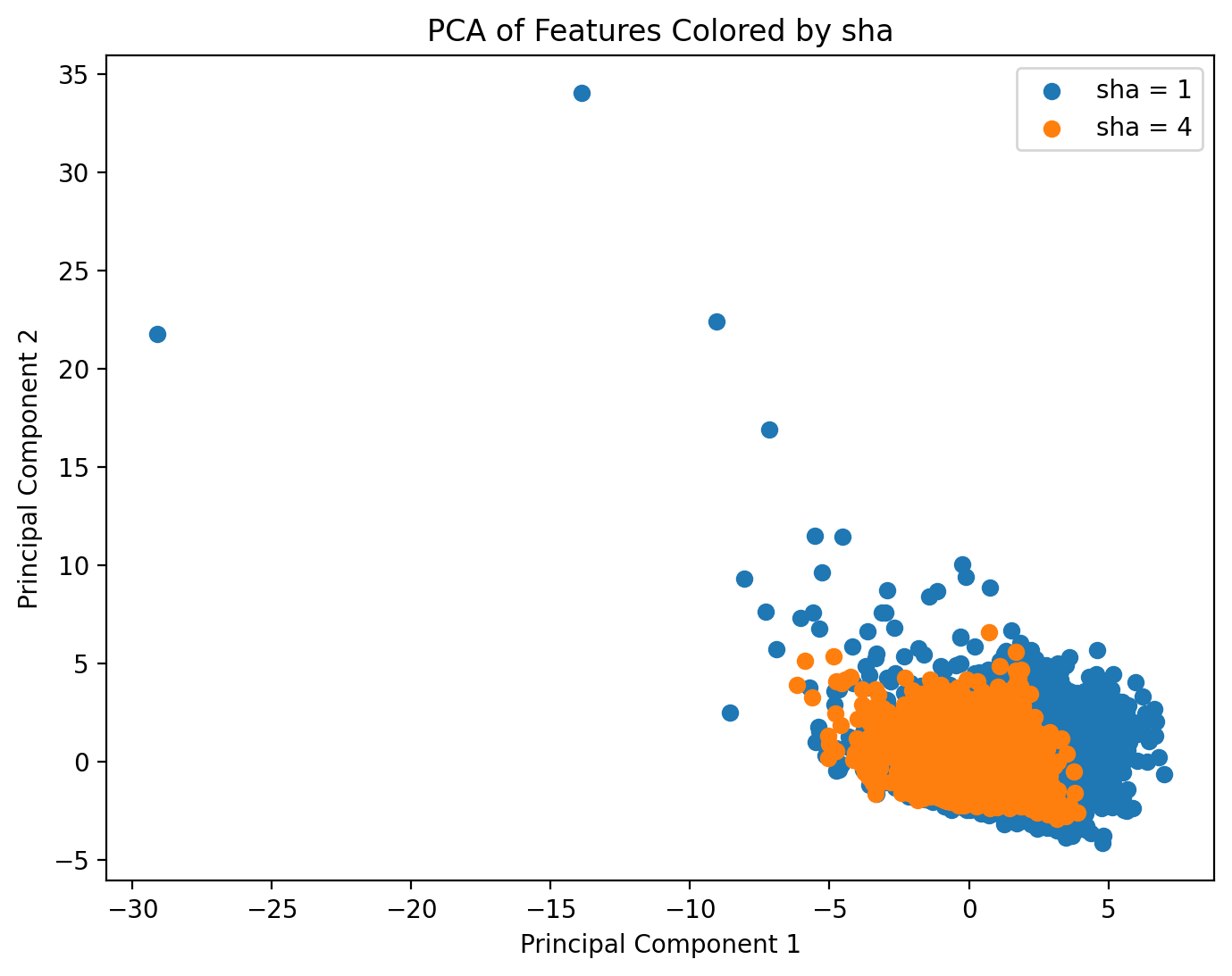}
    \caption{PCA on the dataset used in \Cref{sec:1-vs-4}.}
    \label{fig:pca_1_4}
\end{figure}

In \Cref{exp: 1-vs-4} it was observed that removing the Tamagawa product as a feature from the training set decreased the accuracy of the model, suggesting that this is a predictive feature of the size of $\Sha(E/\Q)$. To further investigate this, we conduct a Principal Component Analysis on this balanced dataset. We plot PC1 and PC2 in \Cref{fig:pca_1_4}, distinguishing between $|\Sha(E/\Q)| = 1$ and $|\Sha(E/\Q)| = 4$.

The large overlap between the two classes in the PCA space shows that the first two principal components are not effective at separating the classes. The loadings for the PCA are given in \Cref{tb:pc_loadings_1_4}. 

\begin{table}[h!]
\centering
\begin{tabular}{lcc}
\hline
\textbf{Feature}         & \textbf{PC1} & \textbf{PC2} \\
\hline
tamagawa\_product & -0.13 &  0.32 \\
rank              &  0.54 &  0.28 \\
conductor         & -0.08 &  0.55 \\
special\_value    &  0.19 &  0.67 \\
torsion           & -0.51 &  0.08 \\
regulator         & -0.19 &  0.14 \\
real\_period      &  0.59 & -0.22 \\
\hline
\end{tabular}
\caption{PCA loadings for the first two principal components}
\label{tb:pc_loadings_1_4}
\end{table}

We see that for PC1, the real period and rank contribute positively, while torsion contributes negatively; these are the dominant features for PC1, and suggest an underlying relationship between them. For PC2, the dominant features are conductor and special value.

\begin{figure}[htbp] 
    \centering
    \includegraphics[width=0.5\linewidth]{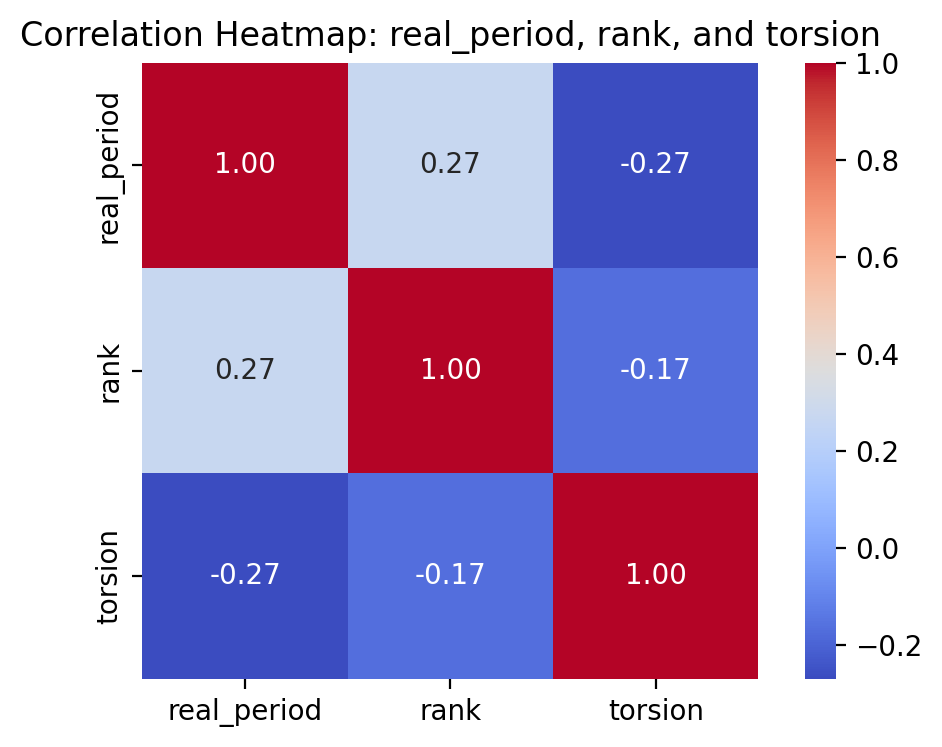}
    \caption{Correlation heatmap on the entire LMFDB dataset between the three features $\Omega$, $r$ and $|E(\Q)_{tors}|$.}
    \label{fig:heatmap_torsion_rank_real_period}
\end{figure}

To investigate any correlation between $\Omega$, $r$ and $|E(\Q)_{tors}|$, we use the entire dataset (not merely the one with positive rank and $\Sha$ size 1 or 4) and start by seeing the correlation coefficients in \Cref{fig:heatmap_torsion_rank_real_period}.

None of these values are particularly high, suggesting that none of these three features are pairwise correlated.

\subsection{PCA for the entire dataset}

\begin{figure}[htbp] 
    \centering
    \includegraphics[width=\linewidth]{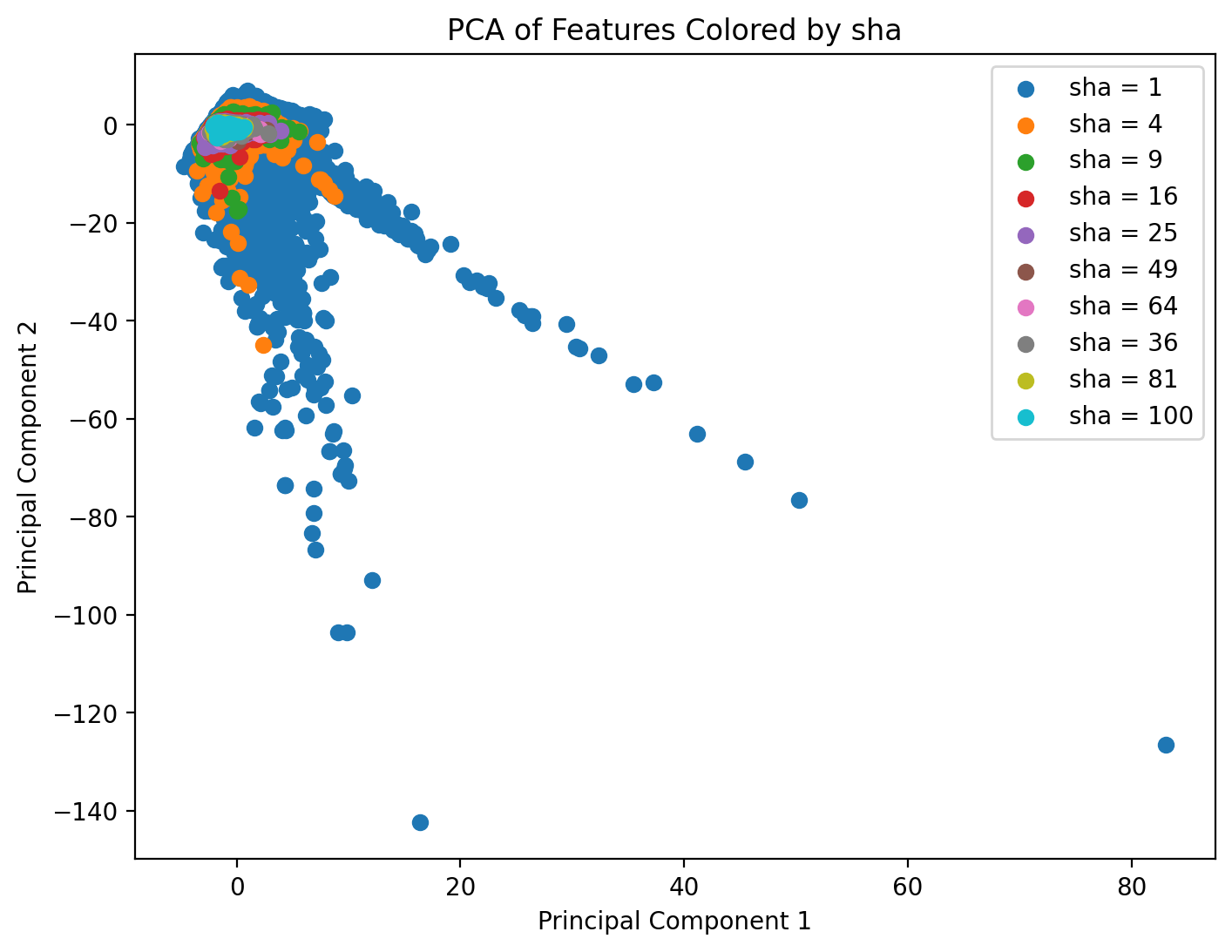}
    \caption{PCA on the dataset with the $10$ largest values for $|\Sha|$.}
    \label{fig:pca_top_10}
\end{figure}

For completeness, we provide a PCA conducted on the dataset with the 10 largest values of $\Sha$ size. This is shown in \Cref{fig:pca_top_10}, with the loadings given in \Cref{tb:pc_loadings_top_10}.

\begin{table}[h!]
\centering
\caption{PCA Loadings for the First Two Principal Components}
\begin{tabular}{lcc}
\hline
\textbf{Feature}         & \textbf{PC1} & \textbf{PC2} \\
\hline
tamagawa\_product &  0.05 & -0.49 \\
rank              &  0.65 &  0.04 \\
conductor         &  0.20 & -0.30 \\
special\_value    &  0.66 & -0.13 \\
torsion           & -0.18 & -0.48 \\
regulator         &  0.17 & -0.27 \\
real\_period      &  0.18 &  0.59 \\
\hline
\end{tabular}
\label{tb:pc_loadings_top_10}
\end{table}

It is curious that there are two distinct `arms' visible in the plot, for which we can find no explanation.

\section{Future work}\label{sec:future}

This work has followed in the footsteps of \cite{alessandretti2023machine} and \cite{he2023curves} and made certain improvements to these works as outlined in \Cref{ssec:main_contributions}. There is much scope for further improvement. In this section we limit ourselves to mentioning three possible future avenues of investigation.

\begin{enumerate}
    \item \textbf{ML for record-breaking $\Sha$}. The largest $\Sha(E)$ known (under BSD) has size $1{,}029{,}212^2$, from the elliptic curve 
    \begin{align*}
        E : y^2 = x^3 &+ 212710514871660026303688x^2 \\
        &+ 11311440784241675303955251133401838632693717904x.
    \end{align*}
    This curve is from \cite{dabrowski2021elliptic}. Can one use a machine learning technique to break this record? (Possibly a model coming from the realm of anomaly detection.) 
    \item \textbf{Neural network regression model}. As explained in \Cref{sec:regression_model}, we were unable to find a suitable neural network architecture that performed well for the regression problem. It may be instructive to consider the models developed by Kazalicki and Vlah \cite{kazalicki2023ranks}, who trained neural networks for learning the rank of elliptic curves. 
    \item \textbf{More advanced data visualisation techniques}. In \Cref{sec:pca} we limited ourselves to PCA, which is one of the simplest approaches to dimensionality reduction. It would be interesting to see how more advanced techniques, perhaps coming from the realm of topological data analysis, could be used to find separation between $\Sha$ classes. One candidate for this could be \emph{IsUMap} \cite{barth2024isumapmanifoldlearningdata}, which was presented at the closing workshop of the Harvard CMSA program on Mathematics and Machine Learning.
\end{enumerate}

\bibliographystyle{alpha}
\bibliography{references.bib}{}
\end{document}